\documentclass[11pt,reqno]{amsart}
\usepackage{graphicx}
\usepackage{amsthm,hyperref}
\newtheorem{theorem}{Theorem}[section]
\newtheorem{lemma}[theorem]{Lemma}
\theoremstyle{definition}
\newtheorem{definition}[theorem]{Definition}
\newtheorem{proposition}{Proposition}

\theoremstyle{remark}
\newtheorem{remark}[theorem]{Remark}
\numberwithin{equation}{section}

\begin{document}
\title[Nonlinear problem by Topological degree]{Nonlinear problem involving the fractional $p(x)-$Laplacian operator by Topological degree}
\author{M. AIT HAMMOU}
\address[M. AIT HAMMOU]{Sidi Mohamed Ben Abdellah university, Laboratory LAMA, Department of Mathematics, Fez, Morocco}
\email[M. AIT HAMMOU]{mustapha.aithammou@usmba.ac.ma}
\subjclass[2010]{Nonlinear elliptic problem, fractional $p(x)-$Laplacian operator, fractional Sobolev spaces with variable exponent, Degree theory.}
\keywords{35J66, 47H30, 46E35, 35R11, 47H11.}

\begin{abstract}
This paper is concerned with the study of the nonlinear problems involving the fractional $p(x)-$Laplacian operator
\begin{equation*}
\left\{\begin{array}{ccc}
(-\Delta_{p(x)})^su(x)+|u(x)|^{q(x)-2}u(x)=\lambda |u(x)|^{r(x)-2}u(x) & \mbox{in}\; \Omega,\\
u=0 &\mbox{in}\; \mathbb{R}^N\setminus\Omega.
\end{array}\right.
\end{equation*}
By means of the Berkovits degree theory, we prove the existence of  nontrivial weak solutions for this problem. The appropriate functional framework for this problems is the fractional Sobolev spaces  with variable exponent.
\end{abstract}
\maketitle

\section{Introduction}
Great attention was paid to the study of elliptic problems introducing the fractional operator during the last years. The study of this type of problems is motivated by by their abilities to model several physical phenomena such as those of phase transition, continuum mechanics and dynamics. Fractional operator is also present in game theory and probability which he provides a simple model to describe stochastic stabilization of L\'evy jump processes (see \cite{A,BRS,BV,C} and the references therein).\par
Let $\Omega$ be a smooth bounded open set in $\mathbb{R}^N$, $s\in (0,1)$ and let $p:\overline{\Omega}\times\overline{\Omega}\rightarrow (1,+\infty)$ be a continuous bounded function. We assume that
\begin{equation}\label{(1)}
 1<p^-=\min_{(x,y)\in\overline{\Omega}\times\overline{\Omega}}p(x,y)\leq p(x,y)\leq p^+=\max_{(x,y)\in\overline{\Omega}\times\overline{\Omega}}p(x,y)<+\infty,
\end{equation}
and $p$ is symmetric i.e.
\begin{equation}\label{sym}
  p(x, y) = p(y, x),\;\;\forall(x, y)\in \overline{\Omega}\times\overline{\Omega}.
\end{equation}
Let denote by :
$$q(x)=p(x,y),\;\; \forall x\in\overline\Omega.$$
Let us consider the fractional $p(x)-$Laplacian operator given by
$$(-\Delta_{p(x)})^su(x)=p.v.\int_\Omega\frac{|u(x)-u(y)|^{p(x,y)-2}(u(x)-u(y))}{|x-y|^{N+sp(x,y)}}\;dy ,\;\;\forall x\in\Omega,$$
where $p.v.$ is a commonly used abbreviation in the principal value sense.\par
In this paper, we are concerned with the study of the the following nonlinear elliptic problem,
\begin{equation}\label{Pr1}\tag{P}
\left\{\begin{array}{ccc}
(-\Delta_{p(x)})^su(x)+|u(x)|^{q(x)-2}u(x)=\lambda |u(x)|^{r(x)-2}u(x) & \mbox{in}\; \Omega,\\
u=0 &\mbox{in}\; \mathbb{R}^N\setminus\Omega.
\end{array}\right.
\end{equation}
where  $r(\cdot)\in C(\bar{\Omega})$ and $\lambda$ is a real parameter. We assume also that
 \begin{equation}\label{(3)}
   1<r^-\leq r(x)\leq r^+<p^-.
 \end{equation}
Note that $(-\Delta_{p(x)})^s$ is the fractional version of well known $p(x)-$Laplacian operator $-\Delta_{p(x)}(u)=-div(|\nabla u|^{p(x)-2}\nabla u)$
for which Alsaedi in \cite{Als} establishes sufficient conditions for the existence of nontrivial weak solutions for a problem similar to \eqref{Pr1} that is the following problem:
\begin{equation*}
  \left\{\begin{array}{ccc}
           -\Delta_{p(x)}u=\lambda |u|^{p(x)-2}u+|u|^{q(x)-2}u& \mbox{ in }\; \Omega, \\
           u=0 &\mbox{ on }\; \partial\Omega.
         \end{array}\right.
\end{equation*}
The proofs combine the Ekeland variational principle, the mountain pass theorem and energy arguments.\par
In \cite{ABS}, the autors study the problem \eqref{Pr1} as an eigenvalue problem using adequate variational techniques, mainly based on Ekeland's variational principle. They establish the existence of a continuous family of eigenvalues lying in a neighborhood at the right of the origin.\par
Using another technical approach: that of the topological degree theory, notably the recent Berkovits degree, we prove in this paper the existence of at least one nontrivial weak solution the problem \eqref{Pr1}.\par
 The paper is divided into four sections. In the second section, we introduce some preliminary results about Lebesgue and fractional Sobolev spaces with variable exponent, some classes of operators and an outline of the recent Berkovits degree. The third section is reserved for some technical lemmas. Finally, in the fourth section we give our main results concerning the weak solutions of the problem (\ref{Pr1}).
\section{Some preliminary results}
\subsection{Lebesgue and fractional Sobolev spaces with variable exponent}
In this subsection, we first recall some useful properties of the variable exponent Lebesgue spaces $L^{p(x)}(\Omega)$. For more details we refer the reader to \cite{FZ1,KR,ZQF} for more details.\par
Denote $$C_+(\bar{\Omega})=\{h\in C(\bar{\Omega})|\inf_{x\in\bar{\Omega}}h(x)>1\}.$$
For any $h\in C_+(\bar{\Omega})$, we define $$h^+:=max\{h(x),x\in\bar{\Omega}\}, h^-:=min\{h(x),x\in\bar{\Omega}\}.$$
For any $p\in C_+(\bar{\Omega})$ we define the variable exponent Lebesgue space
$$L^{p(x)}(\Omega)=\{u;\ u:\Omega \rightarrow \mathbb{R} \mbox{  is measurable and }\int_\Omega| u(x)|^{p(x)}\;dx<+\infty \}.$$
Endowed with {\it Luxemburg norm} $$\|u\|_{p(x)}=\inf\{{\lambda>0}/\rho_{p(\cdot)}(\frac{u}{\lambda})\leq 1\}.$$ where $$\rho_{p(\cdot)}(u)=\int_\Omega|u(x)|^{p(x)}\;dx, \;\;\; \forall u\in L^{p(x)}(\Omega),$$
$(L^{p(x)}(\Omega),\|\cdot\|_{p(x)})$ is a Banach space, separable and reflexive. Its conjugate space is $L^{p'(x)}(\Omega)$ where
$\frac{1}{p(x)}+\frac{1}{p'(x)}=1$ for all $x\in \Omega.$\\
We have also the following result
\begin{proposition}\label{prop1}
 For any $u\in L^{p(x)}(\Omega)$  we have
 \begin{description}
   \item[(i)] $\|u\|_{p(x)}<1(=1;>1)\;\;\;\Leftrightarrow \;\;\;\rho_{p(\cdot)}(u)<1(=1;>1),$
   \item[(ii)] $\|u\|_{p(x)}\geq1\;\;\;\Rightarrow\;\;\;\|u\|_{p(x)}^{p^-}\leq\rho_{p(\cdot)}(u)\leq\|u\|_{p(x)}^{p^+},$
   \item[(iii)] $\|u\|_{p(x)}\leq1\;\;\;\Rightarrow\;\;\;\|u\|_{p(x)}^{p^+}\leq\rho_{p(\cdot)}(u)\leq\|u\|_{p(x)}^{p^-}.$
 \end{description}
\end{proposition}
From this proposition, we can deduce the inequalities
\begin{equation}\label{e1}
  \|u\|_{p(x)}\leq \rho_{p(\cdot)}(u)+1,
\end{equation}
\begin{equation}\label{ineq1}
  \rho_{p(\cdot)}(u)\leq\|u\|_{p(x)}^{p^-}+\|u\|_{p(x)}^{p^+}.
\end{equation}
 If $p, q\in C_+(\overline{\Omega}$ such that $p(x) \leq q(x)$ for any $x\in\bar{\Omega},$ then there exists the continuous embedding $L^{q(x)}(\Omega)\hookrightarrow L^{p(x)}(\Omega).$\par
Next, we present the definition and some results on fractional Sobolev spaces with variable exponent that was introduced in \cite{ABS,BR,KRV}. Let $s$ be a fixed real number such that $0 < s < 1$ and lets the assumptions \eqref{(1)} and \eqref{sym} be satisfied, we define the fractional Sobolev space with variable exponent via the Gagliardo approach as follows:
\begin{eqnarray*}
  W &=&  W^{s,p(x,y)}(\Omega) \\
   &=& \{u\in L^{q(x)}(\Omega):\int_{\Omega\times\Omega}\frac{|u(x)-u(y)|^{p(x,y)}}{\lambda^{p(x,y)}|x-y|^{N+sp(x,y)}}\;dxdy<+\infty , \mbox{ for some } \lambda>0\},
\end{eqnarray*}
where $q(x)=p(x,x)$. We equip the space $W$ with the norm
$$\|u\|_W = \|u\|_{q(x)}+[u]_{s,p(x,y)},$$ where $[\cdot]_{s,p(x,y)}$ is a Gagliardo seminorm with variable exponent, which is defned by
$$[u]_{s,p(x,y)}=\inf\{{\lambda>0}:\int_{\Omega\times\Omega}\frac{|u(x)-u(y)|^{p(x,y)}}{\lambda^{p(x,y)}|x-y|^{N+sp(x,y)}}\;dxdy\leq1 \}.$$
The space $(W,\|\cdot\|_W)$  is a Banach space (see \cite{CLR}), separable and reflexive (see \cite[Lemma 3.1]{BR}).\par
We also define $W_0$ as the subspace of $W$ which is the closure of $C_0^\infty(\Omega)$ with respect to the norm $\|\cdot\|_W$. From \cite[Theorem 2.1 and Remark 2.1]{ABS},
$$\|\cdot\|_{W_0}:=[\cdot]_{s,p(x,y)}$$
is a norm on $W_0$ which is equivalent to the norm $\|\cdot\|_W$, and we have the compact embedding $W_0\hookrightarrow\hookrightarrow L^{q(x)}(\Omega)$. So the space $(W_0,\|\cdot\|_{W_0})$ is a Banach space separable and reflexive.\\
We defne the modular $\rho_{p(\cdot,\cdot)}:W_0\rightarrow\mathbb{R}$ by
$$\rho_{p(\cdot,\cdot)}(u)=\int_{\Omega\times\Omega}\frac{|u(x)-u(y)|^{p(x,y)}}{|x-y|^{N+sp(x,y)}}\;dxdy.$$
The modular $\rho_{p(\cdot,\cdot)}$ checks the following results, which is similar to Proposition \ref{prop1}(see \cite[Lemma 2.1]{ZZ})
\begin{proposition}\label{prop2}
 For any $u\in W_0$ we have
 \begin{description}
   \item[(i)] $\|u\|_{W_0}\geq1\;\;\;\Rightarrow\;\;\;\|u\|_{W_0}^{p^-}\leq\rho_{p(\cdot,\cdot)}(u)\leq\|u\|_{W_0}^{p^+},$
   \item[(ii)] $\|u\|_{W_0}\leq1\;\;\;\Rightarrow\;\;\;\|u\|_{W_0}^{p^+}\leq\rho_{p(\cdot,\cdot)}(u)\leq\|u\|_{W_0}^{p^-}.$
 \end{description}
\end{proposition}
\subsection{Some classes of operators and an outline of Berkovits degree}
 Let $X$ be a real separable reflexive Banach space with dual $X^*$ and with continuous pairing $\langle.\;,\;.\rangle$ and let $\Omega$ be a nonempty subset of $X$.
 \par
Let  $Y$ be a real Banach space. We recall that a mapping $F:\Omega\subset X \rightarrow Y$ is \it bounded, \rm if it takes any bounded set into a bounded set. $F$ is said to be \it demicontinuous, \rm if for any $(u_n)\subset\Omega$, $u_n \rightarrow u$ implies $F(u_n)\rightharpoonup F(u)$. $F$ is said to be \it compact \rm if it is continuous and the image of any bounded set is relatively compact.\par
A mapping $F:\Omega\subset X\rightarrow X^*$ is said to be \it of class $(S_+)$\rm, if for any $(u_n)\subset\Omega$ with $u_n \rightharpoonup u$ and $limsup\langle Fu_n,u_n - u\rangle\leq 0 $, it follows that $u_n \rightarrow u$.\\ $F$ is said to be \it quasimonotone \rm, if for any $(u_n)\subset\Omega$ with $u_n\rightharpoonup u$, it follows that $limsup\langle Fu_n,u_n - u\rangle\geq 0$.\par
For any operator $F:\Omega\subset X\rightarrow X$ and any bounded operator $T:\Omega_1\subset X\rightarrow X^*$ such that $\Omega\subset\Omega_1$, we say that $F$ satisfies condition $(S_+)_T$, if for any $(u_n)\subset\Omega$ with $u_n \rightharpoonup u$, $y_n:=Tu_n\rightharpoonup y$ and
  $limsup\langle Fu_n,y_n-y\rangle\leq 0$, we have $u_n\rightarrow u$. We say that $F$ has the property $(QM)_T$, if for any  $(u_n)\subset\Omega$ with $u_n\rightharpoonup u$, $y_n:=Tu_n\rightharpoonup y$, we have $limsup\langle Fu_n,y-y_n\rangle \geq0$.\par
  Let $\mathcal{O}$ be the collection of all bounded open set in $X$. For any $\Omega\subset X$, we consider the following classes of operators:
\begin{eqnarray*}
\mathcal{F}_1(\Omega)&:=&\{\small{ F:\Omega\rightarrow X^*\mid F \mbox{ is bounded, demicontinuous and satifies condition }(S_+)}\}, \\
\mathcal{F}_{T,B}(\Omega)&:=&\{\small{F:\Omega\rightarrow X\mid F \mbox{ is bounded, demicontinuous and satifies condition }(S_+)_T}\}, \\
\mathcal{F}_T(\Omega)&:=&\{F:\Omega\rightarrow X\mid F \mbox{ is demicontinuous and satifies condition }(S_+)_T\},\\
\mathcal{F}_B(X)&:=&\{F \in \mathcal{F}_{T,B}(\bar{G})\mid G\in\mathcal{O}, \rm T\in \mathcal{F}_1(\bar{G})\}.
\end{eqnarray*}
Here, $T\in \mathcal{F}_1(\bar{G})$ is called an \it essential inner map \rm to $F$.
\begin{lemma}\label{l2.1}\cite[Lemmas 2.2 and 2.4]{Ber}
   Suppose that $T\in \mathcal{F}_1(\bar{G})$ is continuous and\\
   $S:D_S\subset X^*\rightarrow X$ is demicontinuous such that $T(\bar{G})\subset D_S$, where $G$ is a bounded open set in a real reflexive Banach space $X$. Then the following statements are true:
   \begin{description}
     \item[(i)] If $S$ is quasimonotone, then $I+SoT\in \mathcal{F}_T(\bar{G})$, where $I$ denotes the identity operator.
     \item[(ii)] If $S$ is of class $(S_+)$, then $SoT\in \mathcal{F}_T(\bar{G})$
   \end{description}
 \end{lemma}
\begin{definition}
   Let $G$ be a bounded open subset of a real reflexive Banach space $X$, $T\in \mathcal{F}_1(\bar{G})$ be continuous and let $F,S\in \mathcal{F}_T(\bar{G})$. The
   affine homotopy \\$H:[0,1]\times\bar{G}\rightarrow X$ defined by $$H(t,u):=(1-t)Fu+tSu \mbox{ for } (t,u)\in[0,1]\times\bar{G}$$ is called an admissible
   affine homotopy with the common continuous essential inner map $T$.
 \end{definition}
\begin{remark}\cite{Ber}
    The above affine homotopy satisfies condition $(S_+)_T$.
 \end{remark}

 We introduce the topological degree for the class $\mathcal{F}_B(X)$ due to Berkovits \cite{Ber}.
 \begin{theorem}\label{t2.1}
   There exists a unique degree function
    $$d:\{(F,G,h)|G\in \mathcal{O}, T\in \mathcal{F}_1(\bar{G}), F\in \mathcal{F}_{T,B}(\bar{G}),h\notin F(\partial G)\}\rightarrow \mathbb{Z}$$
    that satisfies the following properties
   \begin{enumerate}
               \item (Existence) if $d(F,G,h)\neq 0$ , then the equation $Fu=h$ has a solution in $G$.
               \item (Additivity) Let $F\in \mathcal{F}_{T,B}(\bar{G})$. If $G_1$ and $G_2$ are two disjoint open subsets of $G$ such that
               $h\not\in F(\bar{G}\setminus (G_1\cup G_2))$, then we have $$d(F,G,h)=d(F,G_1,h)+d(F,G_2,h).$$
               \item (Homotopy invariance) If $H: [0,1]\times \bar{G}\rightarrow X$ is a bounded admissible affine homotopy with a common continuous essential inner map and $h:[0,1]\rightarrow X$ is a continuous path in $X$ such that $h(t)\notin H(t,\partial G)$ for all $t\in [0,1]$ ,then the value of $ d(H(t,.),G,h(t))$ is constant for all $t\in[0,1]$.
               \item (Normalization) For any $h\in G$, we have $d(I,G,h)=1.$
    \end{enumerate}
 \end{theorem}
\section{Technical lemmas}
Let $\Omega\subset \mathbb{R}^N$, $N\geq2,$ be a smooth bounded open set, $s\in(0,1)$ and we assume that \eqref{(1)}, \eqref{sym} and \eqref{(3)} holds. In this section, we present two technical lemmas that we will need to study our problem \eqref{Pr1}.\par
Let denote $L:W_0\rightarrow W_0^*,$ the operator associated to the $(-\Delta_{p(x)})^s$ defined by
$$\langle Lu,v\rangle=\int_{\Omega\times\Omega}\frac{|u(x)-u(y)|^{p(x,y)-2}(u(x)-u(y))(v(x)-v(y))}{|x-y|^{N+sp(x,y)}}|\nabla u|^{p(x)-2}\;\; dxdy,$$ for all$u,v\in W_0,$ where $W_0^*$ is the dual space of $W_0$.
\begin{lemma}\cite{BR}\label{t3.1}
  \begin{description}
    \item[(i)] $L$ is bounded and strictly monotone operator,
    \item[(ii)] $L$ is a mapping of type $(S_+)$,
    \item[(iii)] $L$ is a homeomorphism.
  \end{description}
\end{lemma}
\begin{lemma}\label{l4.1}
 The operator $S:W_0\rightarrow W_0^*$ setting by
  $$\langle Su,v\rangle=\int_\Omega|u(x)|^{q(x)-2}u(x)v(x)\; dx-\lambda\int_\Omega|u(x)|^{r(x)-2}u(x)v(x) dx,\;\;\ \forall u,v\in W_0$$ is compact.
\end{lemma}
\proof
  \underline{\bf Step 1}\\
  Let $\phi:W_0\rightarrow L^{q'(x)}(\Omega)$ be the operator defined by
   $$\phi u(x):=|u(x)|^{q(x)-2}u(x) \mbox{ for } u\in W_0 \mbox{ and } x\in \Omega.$$
  It's obvious that $\phi$ is continuous. We prove that $\phi$ is bounded.\\
For each $u\in W_0$, we have by the inequalities \eqref{e1} and \eqref{ineq1} that
\begin{eqnarray*}
  \|\phi u\|_{q'(x)} &\leq& \rho_{q'(\cdot)}(\phi u)+1 \\
   &=& \int_\Omega||u|^{q(x)-1}|^{q'(x)}\;dx+1 \\
   &=& \rho_{q(\cdot)}(u) \\
   &\leq& \|u\|_{q(x)}^{q^-}+\|u\|_{q(x)}^{q^+})+1.
\end{eqnarray*}
By the compact embedding $W_0\hookrightarrow\hookrightarrow L^{q(x)}(\Omega)$, we have
$$\|\phi u\|_{q'(x)}\leq const(\|u\|_{W_0}^{q^-}+\|u\|_{W_0}^{q^+})+1$$
This implies that $\phi$ is bounded on $W_0$.\\
\underline{\bf Step 2}\\
 Let $\psi:W_0\rightarrow L^{q'(x)}(\Omega)$ be an operator defined by
   $$\psi u(x):=-\lambda |u(x)|^{r(x)-2}u(x) \mbox{ for } u\in W_0 \mbox{ and } x\in \Omega.$$
 $\psi$ is also obviously continuous. By an analog proof in the first step, and in virtu of the continuous embedding $L^{q(x)}(\Omega)\hookrightarrow L^{\alpha(x)}(\Omega)$ where $\alpha(x)=(r(x)-1)q'(x)\in C_+(\overline{\Omega})$ with $\alpha(x)\leq q(x)$ , we obtain
$$\|\psi u\|_{q'(x)}\leq const((\|u\|_{W_0}^{\alpha^-}+\|u\|_{W_0}^{\alpha^+})+1$$
This implies that $\phi$ is bounded on $W_0.$\\
\underline{\bf Step 3}\\
Since the embedding $I:W_0\rightarrow L^{q(x)}(\Omega)$ is compact, it is known that the adjoint operator $I^*:L^{q'(x)}(\Omega)\rightarrow W_0^*$ is also compact. Therefore, the compositions\\
$I^*o\phi \mbox{ and } I^*o\psi :W_0\rightarrow W_0^*$ are compact. We conclude that $S=I^*o\phi+I^*o\psi$ is compact.
\endproof
\section{Main Result}
In this section, we study the nonlinear problem (\ref{Pr1}) based on the Berkovits degree theory introduced in subsection 2.2, where $\Omega\subset \mathbb{R}^N$, $N\geq2,$ is a smooth bounded open, $s\in(0,1)$ and under assumptions \eqref{(1)}, \eqref{sym} and \eqref{(3)}.\par
Let $L \mbox{ and } S:W_0\rightarrow W_0^*(\Omega)$ be as in Section 3.
\begin{definition}
  We say that $u\in W_0$ is a weak solution of (\ref{Pr1}) if
  $$\langle Lu,v\rangle+\langle Su,v\rangle=0,\;\;\forall v\in W_0.$$
\end{definition}

\begin{theorem}
  Under assumptions \eqref{(1)}, \eqref{sym} and \eqref{(3)}, the problem (\ref{Pr1}) has a weak solution $u$ in $W_0.$
\end{theorem}
\proof
$u\in W_0$ is a weak solution of \eqref{Pr1} if and only if
  \begin{equation}\label{eq2}
    Lu=-Su.
  \end{equation}
Thanks to the properties of the operator $L$ seen in Lemma \ref{t3.1} and in view of Minty-Browder Theorem \cite[Theorem 26A]{Z}, the inverse operator $T:=L^{-1}:W_0^*\rightarrow W_0$ is bounded, continuous and satisfies condition $(S_+)$. Moreover, note by Lemma \ref{l4.1} that the operator $S$ is bounded, continuous and quasimonotone.\\
Consequently, equation (\ref{eq2}) is equivalent to
\begin{equation}\label{eq3}
  u=Tv \mbox{ and } v+SoTv=0.
\end{equation}
To solve (\ref{eq3}), we will apply the degree theory introduced in section 2. To do this, we first claim that the set
$$B:=\{v\in W_0^*|v+tSoTv=0 \mbox{ for some } t\in[0,1]\}$$ is bounded. Indeed, let $v\in B$. Set $u:=Tv$, then $\|Tv\|_{W_0}=\|u\|_{W_0}.$\\
If $\|u\|_{W_0}\leq 1$, then $\|Tv\|_{W_0}$ is bounded.\\
If $\|u\|_{W_0}>1$, then we get by the implication (i) in Proposition \ref{prop2} and the inequality \eqref{ineq1} the estimate
\begin{eqnarray*}
  \|Tv\|_{W_0}^{p^-} &=& \|u\|_{W_0}^{p-} \\
   &\leq& \rho_{p(\cdot,\cdot)}(u) \\
   &=& \langle Lu,u\rangle  \\
  &=& \langle v,Tv\rangle \\
   &=& -t\langle SoTv,Tv \rangle  \\
   &\leq& t\int_\Omega (|u(x)|^{q(x)}+\lambda|u(x)|^{r(x)})u\;dx \\
  &\leq& const(\|u\|_{q(x)}^{q^-}+\|u\|_{q(x)}^{q^+}+\|u\|_{r(x)}^{r^-}+\|u\|_{r(x)}^{r^+}).
\end{eqnarray*}
From the continuous embedding $L^{q(x)}(\Omega)\hookrightarrow L^{r(x)}(\Omega)$ and the compact embedding $W_0\hookrightarrow\hookrightarrow L^{q(x)}(\Omega)$, we can deduct the estimate
$$\|Tv\|_{W_0}^{p^-} \leq const(\|Tv\|_{W_0}^{q^+}+\|Tv\|_{W_0}^{r^+}).$$
It follows that $\{Tv|v\in B\}$ is bounded.\\
Since the operator $S$ is bounded, it is obvious from (\ref{eq3}) that the set $B$ is bounded in $W_0^*$. Consequently, there exists $R>0$ such
that $$\|v\|_{W_0^*}<R \mbox{ for all } v\in B.$$
This says that
$$v+tSoTv\neq0 \mbox{ for all } v\in\partial B_R(0) \mbox{ and all } t\in[0,1].$$
From Lemma \ref{l2.1} it follows that
$$I+SoT\in \mathcal{F}_T(\overline{B_R(0)})\mbox{ and } I=LoT\in \mathcal{F}_T(\overline{B_R(0)}).$$
Since the operators $I$, $S$ and $T$ are bounded, $I+SoT$ is also bounded. We conclude that
$$I+SoT\in \mathcal{F}_{T,B}(\overline{B_R(0)})\mbox{ and } I\in\mathcal{F}_{T,B}(\overline{B_R(0)}).$$
Consider a homotopy $H:[0,1]\times\overline{B_R(0)}\rightarrow W_0^*$ given by
$$H(t,v):=v+tSoTv \mbox{ for } (t,v)\in[0,1]\times\overline{B_R(0)}.$$
Applying the homotopy invariance and normalization property of the degree $d$ stated in Theorem \ref{t2.1}, we get
$$d(I+SoT,B_R(0),0)=d(I,B_R(0),0)=1,$$
and hence there exists a point $v\in B_R(0)$ such that
$$v+SoTv=0.$$
We conclude that $u=Tv$ is a weak solution of \eqref{Pr1}.
\endproof

{\small
{\em Authors' addresses}:
{\em Mustapha AIT HAMMOU}, Laboratory LAMA, Faculty of sciences Dhar el Mahraz, Sidi Mohammed ben Abdellah university, Fez, Morocco.
 e-mail: \texttt{mustapha.aithammou@\allowbreak usmba.ac.ma}.

}

\end{document}